\documentclass[12pt]{article}
\usepackage{euscript}
\usepackage{amssymb}
\usepackage{amsmath}
\usepackage{graphicx}
\usepackage{epsfig}
\usepackage{setspace}
\input psfig

\input psfig
        \newcommand\beqn{\begin{eqnarray*}}
        \newcommand\eeqn{\end{eqnarray*}}
        \newcommand\beqa{\begin{eqnarray}}
        \newcommand\eeqa{\end{eqnarray}}
        \newcommand\br{\begin{array}}
        \newcommand\er{\end{array}}
        \newcommand\bes{\begin{equation*}}
        \newcommand\ees{\end{equation*}}
        \newcommand\be{\begin{eqnarray}}
        \newcommand\ee{\end{eqnarray}}
        \newcommand\efig{\end{figure}}
        \newcommand\btab{\begin{table}}
        \newcommand\etab{\end{table}}
        \newcommand\btabl{\begin{tabular}}
        \newcommand\etabl{\end{tabular}}
        \newcommand\bm{\begin{math}}
        \newcommand\ema{\end{math}}
        \newcommand\bi{\begin{itemize}}
        \newcommand\ei{\end{itemize}}
        \newcommand\bd{\begin{description}}
        \newcommand\ed{\end{description}}
        
        \newcommand{\bpr}{\begin{proof}}
        \newcommand{\epr}{\end{proof}}
        \newtheorem{lemma}{\bf Lemma}
        \newtheorem{theorem}{\bf Theorem}

         \newcommand{ \vect }[2]{
                    \left( \begin{array}{c}
                           #1 \\
                           #2
                            \end{array} \right)}
         \newcommand{ \vectt }[3]{
                    \left( \begin{array}{c}
                           #1 \\
                           #2 \\
                           #3
                            \end{array} \right)}
         
        \newcommand\eio{e^{i \theta}}
        
        \newcommand\Eu{\EuScript}
        \newcommand\mrm{\mathrm}

        \newcommand\cm{\CMcal}
        \newcommand\EuB{\Eu B(\ell^2, \ell^2)}
\newcommand\EuBt{\Eu B(\ell^2, \; \ell^2 \times \ell^2)}
        \newcommand\EuBc{\Eu B_c(\ell^2, \ell^2)}
\newcommand\EuBct{\Eu B_c(\ell^2, \ell^2 \times \ell^2)}
\date{}
\begin{document}
\title{Optimal Disturbance Rejection and Robustness for Infinite Dimensional LTV Systems}
\author{Seddik M. Djouadi
{\thanks{S.M. Djouadi is with the
Electrical Engineering \& Computer Science Department, University of Tennessee, Knoxville,
TN 37996-2100. {\tt\small djouadi@eecs.utk.edu}}}
}
\maketitle
\begin{abstract}
In this paper, we consider the optimal disturbance rejection problem for possibly infinite dimensional
linear time-varying (LTV) systems using a framework based on operator algebras of classes of bounded linear
operators. This approach does not assume any state space representation and views LTV systems as causal operators.
After reducing the problem to a shortest distance minimization in a space of bounded
linear operators, duality theory is applied to show existence of optimal solutions, which satisfy a ``time-varying''
allpass or flatness condition. Under mild assumptions the optimal TV controller is shown to be essentially unique.
Next, the concept of M-ideals of operators is used to show that the computation of time-varying (TV) controllers reduces
to a search over compact TV Youla parameters. This involves the norm of a TV compact Hankel operator defined on
the space of causal trace-class 2 operators and its maximal
vectors. Moreover, an operator identity to compute the optimal TV Youla parameter is provided. These results are
generalized to the mixed sensitivity problem for TV systems as well, where it is shown that the optimum is equal to
the operator induced of a TV mixed Hankel-Toeplitz. The final outcome of the approach developed here is that it leads to two tractable
finite dimensional convex optimizations producing estimates to the optimum within desired tolerances, and a method to
compute optimal time-varying controllers.
\end{abstract}
\section*{Mathematical Preliminaries  and Notation}
\begin{itemize}
\item  $\Eu {B} (E, F)$ denotes the
space of bounded linear operators from a Banach space $E$ to a Banach
space $F$, endowed with the operator norm
\beqn
\| A\| :=\sup_{x\in E, \; \|x\| \leq 1} \|Ax\|, \;\; A\; \in\;\Eu {B} (E, F)
\eeqn
\item $\ell^2$ denotes the usual Hilbert space of square summable sequences
with the standard standard norm
\beqn
\| x \|_2^2 := \sum_{j=0}^\infty |x_j|^2,
\;\; x := \bigl(x_0, x_1, x_2, \cdots \bigl) \in \ell^2
\eeqn
and inner product
\beqn
<x, y> = \sum_i \bar{y}_i x_i,\;\;\; y= \bigl(y_0, y_1, y_2, \cdots \bigl) \in \ell^2
\eeqn
%
\item $P_{k}$ the usual truncation operator for some integer $k$, which sets all
outputs after time $k$ to zero, i.e.,
\beqn
P_{k} x = \bigl(x_0, x_1, x_2, \cdots, x_k, 0, 0, \cdots \bigl)
\eeqn
\item An operator $A \in \Eu {B} (\ell^2, \ell^2)$ is said to be causal if
it satisfies the operator equation:
\beqn
P_{k }A P_{k } = P_{k } A, \; \forall k \;\; {\rm positive\;\;integers}
\eeqn
and strictly causal if it satisfies
\beqa
P_{k+1 }A P_{k } = P_{k+1 } A,, \; \forall k \;\; {\rm positive\;\;integers}
\label{2}
\eeqa
\item A partial isometry is an operator $A$ on a Hilbert space which preserves the norm, i.e., $\|Au\| = \|u\|$ for all vectors $u$ in the orthogonal complement of the null space of $A$ \cite{douglas}.
\item A bounded linear functional $f$ on a Banach space $E$ is a real or complex valued linear map such that there exists $c \geq 0$ with $| <f, x>| \leq c \| x\|_{E}$, $\forall x \in E$, where $<f, x>$ denotes the image of $x$ under $f$, and $\| \cdot\|_E$ is the norm on $E$. The dual pace of $E$, denoted $E^\star$, is the space of bounded linear functionals on $E$ under the norm $\| f\|_{E^\star} :=\sup_{x\in E,\; \|x\|_{E}} |<f,x>|$. Note that $E^\star$ is itself a Banach space \cite{douglas}.
\item A sequence $\{x_n\}$ in $E$ is said to converge weakly to $x\in E$ if $|<f, x_n> - <f, x> |\longrightarrow 0$, $\forall f \in E^\star$. A sequence $\{f_n\} \subset E^\star$ is said to converge in the weak$^\star$ topology to $f \in E^\star$ if $|<f_n, x> - <f, x> | \longrightarrow 0$, $\forall x\in E$.
\item A sequence of operators $\{T_n \} \subset \Eu B (E, F)$ converges to an operator $T \in \Eu B (E, F)$ in the operator topology if $\| T_n - T\| \longrightarrow 0$. It is said to converge in the strong operator topology (SOT) if $\| T_n x -T x\|_{F} \longrightarrow 0$, $\forall x \in E$. $\cal K$ denotes the space of compact operators in $\Eu B (E, F)$. An operator $T \in \Eu B (E, F)$ is compact if there exists a sequence of finite rank (i.e., with finite dimensional ranges) operators $\{T_n\} \subset \Eu B (E, F)$, such that, $\| T- T_n\| \longrightarrow 0$ as $n \longrightarrow \infty$ \cite{douglas}.
\end{itemize}
The subscript ``$_c$'' denotes the restriction of a subspace of operators to its intersection
with causal (see \cite{saeks,feintuch} for the definition) operators.
``$^\star$'' stands for the adjoint of an operator or the dual space of a Banach space
depending on the context. $\ominus$ denotes the orthogonal complement. The symbol "$\bigwedge$" denotes the closed linear span, and "$\bigwedge$" the intersection. ``$\rightharpoonup^\star$'' denotes convergence in the weak$^\star$ topology \cite{douglas,luen}. Right and left hand sides are abbreviated as RHS and LHS,
respectively.
\section{Introduction} \label{sec1}
There have been numerous attempts in the literature to generalize ideas about
robust control theory \cite{dft,zhou} to time-varying (TV) systems
(for e.g. \cite{fein1,fein2,feintuch,ravi,Ichikawa,georgiou,dale,lall,peters} and references therein).
In \cite{fein1,fein2} and more recently \cite{feintuch} the authors studied the optimal weighted
sensitivity minimization problem, the two-block problem, and the model-matching problem for LTV systems using
inner-outer factorization for positive operators. They obtained abstract solutions involving the
computation of norms of certain operators, which is difficult since these are infinite
dimensional problems. Moreover, no indication on how to compute optimal linear time-varying (LTV) controllers is
provided. The operator theoretic methods in \cite{fein1,fein2,feintuch} are difficult
to implement and do not provide algorithms which even compute approximately TV controllers. In \cite{tadmor} a state space
extension of the Nehari Theorem to a time-varying system-theoretic setting is developed with parametrization of suboptimal solutions.
In \cite{lall} the authors rely on state space techniques which lead to algorithms based on
infinite dimensional operator inequalities. These methods lead to suboptimal controllers but are
difficult to solve and are restricted to finite dimensional systems. Moreover, they do not allow the degree
of suboptimality to be estimated. An extension of these results to uncertain systems is reported
in \cite{pirie} relying on uniform stability concepts.
\\
In \cite{kha2} both the sensitivity minimization problem in the presence
of plant uncertainty, and robust stability for LTV systems in the $\ell^\infty$ induced norm
is considered. However, their methods could not be extended to the case of systems operating on
finite energy signals.
\\
We believe that it is important to point to the fact that the lack of calculation of
the degree of suboptimality in these methods, can result in arbitrary poor overall
closed-loop performance. This fact is explored in \cite{zo}, where the authors, although in the linear time-invariant (LTI) case, studied a "two-arc"
counter example that shows that in the limit suboptimal controllers can lead to arbitrary poor performance in the
presence of plant uncertainty. It seems that this counter example has not yet received enough attention.
 \\
Analysis of time-varying control strategies for optimal disturbance rejection for known time-invariant plants has been studied in \cite{shda,chap}.
A robust version of these problems was considered in \cite{jeff,kha1} in different induced
norm topologies. All these references showed that for {\it time-invariant nominal} plants and weighting
functions, time-varying control laws offer no advantage over time-invariant ones.
\\
In the book \cite{Foias}, the authors consider various interpolation problems in the stationary (LTI) and nonstationary (LTV) cases. In particular, Nevanlinna-Pick, Hermite-Fejer, Sarason, and Nudelman interpolations and Nehari extension problems are studied.  In the TV case the unit disk is replaced by the set of diagonal matrices such as the weighted shift with spectral radius less than 1. The interpolation condition is replaced by an identity involving diagonal matrices. Two methods are used to solve these problems, a method based on the reduction to the time-invariant interpolation case but with operator valued functions, and another method based on a TV commutant lifting theorem developed by the authors of the book.  The solutions for these problems are provided in terms of infinite dimensional operator identities.
\\
\\
In this paper, we are interested in optimal disturbance rejection for ({\it possibly infinite-dimensional},
i.e., systems with an infinite number of states) LTV systems. These systems have been used as models
in computational linear algebra and in a variety of computational and communication networks as described
in the excellent book \cite{dewilde}. This allows variable number of states which is predominant in networks which can switch on or off
certain parts of the system \cite{dewilde}, and infinite number of states as in distributed parameter systems. We study questions such as
existence and uniqueness of optimal LTV controllers under specific conditions.
What is the corresponding notion for the Hankel operator which solves the optimal $H^\infty$ problem?
What does allpass mean for LTV systems? Is the optimal Youla parameter compact and if so under what conditions?  Is it possible to compute optimal LTV controllers within desired accuracy? Note that compactness
is important since it allows approximations by finite rank operators (matrices), and therefore makes the problem amenable to numerical computation.
\\
\\
Using inner-outer factorizations as defined in \cite{dav,feintuch} with respect of the nest algebra
of lower triangular (causal) bounded linear operators defined on $\ell^2$ we show that the
problem reduces to a distance minimization between a special operator and the nest algebra.
The inner-outer factorization used here holds under weaker assumptions than
\cite{fein1,fein2}, and in fact, as pointed in (\cite{dav} p. 180), is different from the
factorization for positive operators used there. Duality structure and predual formulation of the problem
showing existence of optimal LTV controllers is provided. The optimum is shown to satisfy
a ``TV'' allpass condition quantified in the form of a partial isometry of an operator, therefore, generalizing
the flatness or allpass condition \cite{zames_francis,fran,bruce}. The optimal controller and the corresponding dual operator
are shown to be essentially unique under mild assumptions.
\\
With the use of M-ideals of operators, it is shown that the
computation of time-varying (TV) controllers reduces to a search over compact TV Youla
parameters. Furthermore, the optimum is shown to be equal to the norm of a compact
time-varying Hankel operator defined on the space of causal Hilbert-Schmidt operators. The latter is a ``natural'' analogous to the Hankel operator used in the LTI case.  An operator equation
to compute the optimal TV Youla parameter is also derived. The results obtained here
lead to a pair of dual finite dimensional convex optimizations which approach the
real optimal disturbance rejection performance from both directions not only producing
estimates within desired tolerances, but allowing the computation of optimal
time-varying controllers. The numerical computation involve solving a semi-definite programming problem, and a search
over lower triangular matrices. \\
The results are generalized to the mixed sensitivity problem for TV systems
as well, where it is shown that the optimum is equal to the operator induced of a TV mixed Hankel-Toeplitz
operator generalizing analogous results known to hold in the LTI case \cite{ozbay,GS2,fran}.
\\
Our approach is purely input-output and does not use any state space realization, therefore
the results derived here apply to {\it infinite dimensional LTV systems}, i.e., TV systems with an
infinite number of state variables. Although the theory is
developed for causal stable system, it can be extended in a straightforward fashion to the
unstable case using coprime factorization techniques for LTV systems discussed in
\cite{fein2,feintuch}.
The framework developed can also be applied to other performance indexes, such as the optimal TV robust disturbance attenuation
problem considered in \cite{zo,dj1,dj2}.  For continuous LTV systems subject to time-varying unstructured uncertainty the problem was
considered in \cite{msd8}, where it is shown that for causal  LTV systems, it is equivalent to finding the smallest fixed point of a 'two-disc' type optimization problem under TV feedback control laws. For the discrete version of the same problem, the duality structure of the problem and a solution based on a new bilinear map was provided in terms of an infinite dimensional identity in \cite{djouadiSIAM}. A related problem involving the computation of
the gap metric for LTV systems has been considered in \cite{djouadiSCL}. In particular, the computation of the gap metric for TV systems using a time-varying generalization of normalized coprime factorization as an iterative scheme involving the norm of an operator with a TV Hankel plus Toeplitz structure is provided. The duality structure is characterized by computing the preannihilator and exploiting the particular structure of the problem.
\\
Part of the results presented here were announced in \cite{cdc04_tv,acc07} without proofs.
\\
\\
The rest of the paper is organized as follows. Section \ref{formulation} contains the problem formulation.
The duality structure of the problem is worked out in section \ref{duality}. In section \ref{allpass}, the optimal solution
is shown to satisfy a TV allpass property. Section \ref{mideal} shows that under mild assumptions the TV optimal $Q$ parameter is compact. Section
\ref{uniqueness} discusses the uniqueness of the optimal solution. In section \ref{hankel}, a solution based on a TV
Hankel operator is derived together with an extremal identity for the optimum. Section \ref{numerical} presents a numerical
solution based on duality theory. A Generalization to the TV mixed sensitivity problem is carried out in section \ref{mixed}. Section
\ref{conclusion} contains concluding remarks.
\section{Problem formulation}
\label{formulation}
\begin{figure}[h]
  \begin{center}
    \leavevmode
\psfig{figure=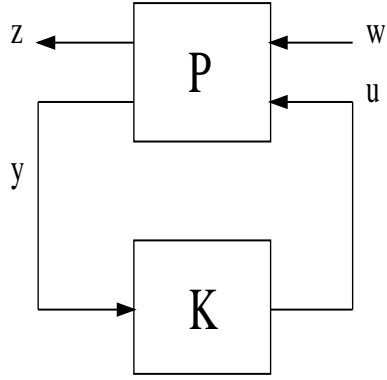,width=5cm,height=5cm,angle=-90}
  \end{center}
\caption{Block Diagram for Disturbance Attenuation} \label{fig1}
\end{figure}

In this paper we consider the problem of optimizing performance for causal
linear time varying systems. The standard block diagram for the optimal
disturbance attenuation problem that is considered here is represented in Fig.
\ref{fig1}, where $u$ represents the control inputs, $y$ the measured outputs,
$z$ is the controlled output, $w$ the exogenous perturbations. $P$ denotes a causal
stable linear time varying plant, and $K$ denotes a time varying controller.
The closed-loop transmission from $w$ to $z$ is denoted by $T_{zw}$. Using the
standard Youla parametrization of all stabilizing controllers 
the closed loop operator $T_{zw}$ can be written as 
\cite{feintuch},
\be
T_{zw} = T_1 - T_2 Q T_3
\ee
%
where $T_1$, $T_2$ and $T_3$ are stable causal time-varying operators, that is,
$T_1$, $T_2$ and $T_3 \in B_c(\ell^2, \ell^2)$. In this paper we assume
without loss of generality that $P$ is stable, the Youla parameter
$Q := K(I+PK)^{-1}$ is then an operator belonging to $B_c(\ell^2, \ell^2)$,
and is related univocally to the controller $K$ \cite{saeks}. Note that $Q$
is allowed to be {\it time-varying}. If $P$ is unstable it suffices to use the
coprime factorization techniques in \cite{dale,feintuch} which lead
to similar results. The magnitude of the signals $w$ and $z$ is measured in the
$\ell^2$-norm. Two problems are considered here optimal disturbance rejection which,
and the mixed sensitivity problem for LTV systems which includes a robustness problem
in the gap metric studied in \cite{feintuch,georgiou}. Note that for the latter problem
$P$ is assumed to be unstable and we have to use coprime factorizations.
The performance index can be written in the following form
\beqa
\mu &:=& \inf \left\{\| T_{zw}\|: \; K \;\; {\rm being \; robustly \;stabilizing\;
linear \; time-varying \; controller} \right\}
\nonumber \\
&=& \inf_{Q \in \Eu{B}_c(\ell^2, \ell^2)} \| T_1 - T_2 Q T_3 \|
\label{4}
\eeqa
The performance index (\ref{4}) will be transformed into a distance minimization between
a certain operator and a subspace to be specified shortly. To this end, define a nest
$\cal{N}$ as a family of closed subspaces of the Hilbert space $\ell^2$ containing
$\{ 0 \}$ and $\ell^2$ which is closed under intersection and closed span.
Let $Q_n := I-P_n,\;\; {\rm for}\;\; n=-1,0,1, \cdots$, where $P_{-1} := 0$
and  $P_\infty := I$. Then $\widetilde{Q}_n$ is a projection, and we can associate to it the
following nest $\mathcal{N} := \{ \widetilde{Q}_{n} {\ell}^2, \; n=-1, 0, 1, \cdots \}$.
The triangular or nest algebra ${\cal{T}} ({\cal{N}})$ is the set of
all operators $T$ such that $TN \subseteq N$ for every element $N$ in ${\cal{N}}$.
That is
\beqa
{\cal{T}} ({\cal{N}}) &=& \{ A \in {\Eu{B}}(\ell^2, \ell^2):\;  P_n A (I-P_n) = 0,
\; \forall \; n\} 
\label{5a}
\eeqa
Note that the Banach space $\Eu{B}_c(\ell^2, \ell^2)$ is identical to the nest algebra
${\cal{T}} ({\cal{N}})$. 
Define
\beqa
N^{-} = \bigvee \{N^\prime \in {\cal{N}} : N^\prime < N \}, \;\;\;\;
N^{+} = \bigwedge \{N^\prime \in {\cal{N}} : N^\prime > N \}
\label{7} 
\eeqa
where $N^\prime < N$ means $N^\prime \subset N$, and $N^\prime > N$ means
$N^\prime \supset N$. The subspaces $N \ominus N^{-}$ are called the atoms of
${\cal{N}}$. Since in our case the atoms of ${\cal{N}}$ span $\ell^2$, then
${\cal{N}}$ is said to be atomic \cite{dav}.
\\
Following \cite{dav} we introduce inner-outer factorizations for the
operators in ${\cal{T}} ({\cal{N}} )$ as follows:
\\ An operator $A$ in ${\cal{T}} ({\cal{N}} )$ is called {\it outer} if the range projection, denoted
$P(R_A)$, $R_A$ being the range of $A$ and $P$ the orthogonal projection onto $R_A$,
commutes with $\mathcal N$ and $AN$ is dense in the intersection $N \cap R_A$ for every $N \in \mathcal N$.
A partial isometry $U$ is called {\it inner} in ${\cal{T}}({\cal{N}} )$ if $U^\star U$
commutes with ${\mathcal{N}}$ \cite{arveson,dav,feintuch}. In our case, $A \in \mathcal
T (\mathcal N) = \Eu{B}_c(\ell^2, \ell^2)$ is outer if $P$ commutes with each $\widetilde{Q}_n$
and $A\widetilde{Q}_n \ell^2$ is dense in $\widetilde{Q}_n \ell^2 \cap A \ell^2$. $U\in \Eu{B}_c(\ell^2, \ell^2)$
is inner if $U$ is a partial isometry and $U^\star U$ commutes with every $\widetilde{Q}_n$. Applying
these notions to the time-varying operator $T_2 \in \Eu{B}_c (\ell^2,\; \ell^2)$,
we get $T_2 = T_{2i}T_{2o}$, where $T_{2i}$ and $T_{2o}$ are inner outer operators in
$\Eu{B}_c (\ell^2, \ell^2)$, respectively. Similarly, outer-inner factorization can
be defined and the operator $T_3= T_{3o}T_{3i}$, where $T_{3i}\in \Eu{B}_c(\ell^2, \ell^2)$ is inner 
$T_{3o}\in \Eu{B}_c(\ell^2, \ell^2)$ is outer.
The performance index $\mu$ in (\ref{4}) can then be written as
\beqa
\mu =\inf_{Q \in \Eu{B}_c(\ell^2, \ell^2)} \| T_1 - T_{2i}T_{2o} Q T_{3o}T_{3i} \|
\label{8}
\eeqa
Following classical robust control theory \cite{fran,bruce,zhou} assume: {\bf (A1)} that $T_{2o}$ and $T_{3o}$ are invertible both in $\Eu{B}_c(\ell^2, \ell^2)$. \\
Assumption (A1) can be somewhat relaxed by assuming instead that the outer operators $T_{2o}$ and $T_{3o}$
are bounded below (see Lemma \cite{arveson} p. 220).
\\
Assumption (A1) guarantees that the map $
Q \longrightarrow T_{2o} \Eu{B}_c(\ell^2, \ell^2) T_{3co}$ is bijective.
In the time-invariant case this assumption means essentially that the
outer factor of the plant $P$ is invertible \cite{bruce}. Under this assumption
$T_{2i}$ becomes an isometry and $T_{3i}$ a co-isometry in which case $
T_{2i}^\star T_{2i} = I$ and $T_{3i} T_{3i}^\star = I$. By ''absorbing'' the operators
$T_{2o}$ and $T_{3co}$ into the ''free'' operator $Q$, expression (\ref{8}) is then equivalent to
\beqa
\mu =\inf_{Q \in \Eu{B}_c(\ell^2, \ell^2)} \| T_{2i}^\star T_1 T_{3ci}^\star - Q  \|
\label{9}
\eeqa
Expression (\ref{9}) is the distance from the operator $T_{2i}^\star T_1 T_{3ci}^\star
\in \Eu{B}(\ell^2, \ell^2)$ to the nest algebra $\Eu{B}_c(\ell^2, \ell^2)$.
In the next section we study the distance minimization problem (\ref{9})
in the context of M-ideals and the operator algebra setting discussed above.
\section{Duality}
\label{duality}
Let $X$ ba a Banach space and $X^\star$ its dual space, i.e., the space of bounded linear
functionals defined on $X$ . For a subset $J$ of $X$, the annihilator of $J$ in
$X^\star$ is denoted $J^\perp$ and is defined by \cite{luen}, $
J^\perp := \{ \Phi \in X^\star :\; \Phi(f) =0, \; f \in J \}$, i.e., $J^\perp$ is the set of bounded linear
functionals on $X$ which vanish on $J$.
Similarly, if $K$ is a subset of $X^\star$ then the preannihilator of $K$ in X is denoted
$^\perp K$, and is defined by $^\perp K := \{ x \in X :\;  \Phi(x) =0, \; \Phi \in K \}$.
The existence of a preannihilator implies that the following identity holds \cite{luen}
\beqa
\min_{y \in K} \| x - y\| = \sup_{k \in ^\perp K, \; \|k\| \leq 1} | <x, \; k>|
\label{15a}
\eeqa
where $<\cdot,\; \cdot>$ denotes the duality product. Let us apply these results to the
problem given in (\ref{9}) by putting
\beqa
X = \Eu{B}(\ell^2,\; \ell^2), \;\; x = T_{2i}^\star T_1 T_{3ci}^\star \in \Eu{B}(\ell^2,\; \ell^2),\;\; J = \Eu{B}_c (\ell^2, \ell^2)
\eeqa
Introduce the class of compact operators on $\ell^2$ called the trace-class or Schatten
1-class, denoted ${\cal C}_1$, such that if $\{\xi_i\}$ is a basis of $\ell^2$, and $T \in {\cal C}_1$, then $\sum_i <T\xi_i, \xi_i> < \infty$. Moreover,
the sum $\sum_i <T\xi_i, \xi_i>$ is independent of the choice of the basis \cite{schatten,dav}. The operator $T^\star T$, $T^\star$ being the adjoint operator of $T$, is positive (i.e., $<Tx, x> >0$, $x\in \ell^2$ and $x \neq 0$) and compact and has a unique positive square root, denoted $(T^\star T)^{\frac{1}{2}}$.  The trace-class norm is defined as \cite{schatten,dav},
\beqa
\| T \|_1 \; :=\; tr (T^\star T)^{\frac{1}{2}} = \sum_i<(T^\star T)^{\frac{1}{2}} \xi_i, \xi_i>
\eeqa
where $tr$ denotes the Trace.
\\
We identify $\Eu{B}(\ell^2, \; \ell^2)$ with the dual space of ${\cal C}_1$ ,
${\cal C}_1^\star$, under trace duality \cite{schatten}, that is, every operator $A$
in $\Eu{B}(\ell^2, \; \ell^2)$ induces a bounded linear functional on ${\cal {C}}_1$
as follows: $\Phi_A \in {\cal C}_1^\star$ is defined by $\Phi_A(T) = tr(AT)$, and we write $\Eu{B}(\ell^2, \; \ell^2) \simeq {\cal C}_1^\star$ to express
that $\Eu{B}(\ell^2, \; \ell^2)$ is isometrically isomorphic to ${\cal C}_1^\star$.
\\
Every trace-class operator $T$ in turn induces a bounded linear functional on
$\Eu{B}(\ell^2, \; \ell^2)$, namely, $\Phi_T(A) = tr (AT)$ for all $A$ in $\Eu{B}(\ell^2, \; \ell^2)$.
\\
To compute the preannihilator of $\Eu{B}_c(\ell^2, \; \ell^2)$ define the subspace $S$ of ${\cal C}_1$ by
\beqa
S := \{ T \in {\cal C}_1:\; (I-\widetilde{Q}_n) T \widetilde{Q}_{n+1}=0, \;\; {\rm for \; all\;} n \}
\label{17}
\eeqa
In the following Lemma we show that $S$ is the preannihilator of $\Eu{B}_c(\ell^2, \; \ell^2)$.
\begin{lemma} \label{l1}
The preannihilator of $\Eu{B}_c(\ell^2, \; \ell^2)$ in ${\cal C}_1$,
$^\perp \Eu{B}_c(\ell^2, \; \ell^2)$, is isometrically isomorphic to $S$.
\end{lemma}
{\bf Proof.}
By Lemma 16.2 in \cite{dav} the preannihilator of ${\cal{T}}({\cal N})$ is given by $\Phi_T \in
^\perp \Eu{B}_c(\ell^2,\; \ell^2)$ if and only $T$ belongs to the subspace
\beqn
\{T \in {\cal C}_1 \;:\; P(N^-)^\perp T P(N) = 0\;\; {\rm for \; all}\; N \; {\rm in}\in {\cal N}\}
\eeqn
where $P(N)$ denotes the orthogonal projection on $N$, likewise for $P(N^-)$.
$P(N^-)^\perp$ the complementary projection of $P(N^-)$, that is, $P(N^-)^\perp = I-P(N^-)$.
\\
In our case ${\cal N}$ is atomic, and for any $N \in {\cal N}$ there exists $n$ such that
$N = \widetilde{Q}_{n+1} \ell^2$, i.e., $P(N) = \widetilde{Q}_{n+1}$. The immediate predecessor of $N$, $N^-$, is then
given by $\widetilde{Q}_n \ell^2$, i.e., $N^- =\widetilde{Q}_n \ell^2$. The orthogonal complement of $N^-$ is then
$\ell^2 \ominus \widetilde{Q}_n \ell^2$. So $P(N^-)^\perp = I- \widetilde{Q}_n$. Therefore,
\beqn
S = \{T \in {\cal C}_1 \;:\; (I-\widetilde{Q}_n) T \widetilde{Q}_{n+1} = 0, \; \forall\; n \}
\eeqn
is isometrically isomorphic to the preannihilator of ${\Eu B_c}(\ell^2,\; \ell^2)$, and the Lemma is proved. $\blacksquare$
\\
\\
The existence of a predual ${\cal C}_1$ and a preannihilator $S$ implies the
following Theorem which is a consequence of Theorem 2 in \cite{luen} (Chapter 5.8).
\begin{theorem} \label{t1}
Under assumption (A1) there exists at least one optimal $Q_o$ in $\Eu{B}_c (\ell^2, \ell^2)$ achieving optimal
performance $\mu$ in (\ref{9}), moreover the following identities hold
\beqa
\mu =\inf_{Q \in \Eu{B}_c(\ell^2, \ell^2)} \| T_{2i}^\star T_1 T_{3ci}^\star - Q  \|
= \|T_{2i}^\star T_1 T_{3ci}^\star - Q_o  \|
= \sup_{T \in S, \; \| T\|_1 \leq 1} |tr(TT_{2i}^\star T_1 T_{3ci}^\star)|
\label{16}
\eeqa
\end{theorem}
Theorem \ref{t1} not only shows the existence of an optimal LTV controller, but plays an important role in its computation by reducing the problem to primal and dual of finite dimensional convex optimizations. Under a certain condition in the next section it is shown that the supremum in (\ref{16}) is achieved.
Theorem \ref{t1} also leads to a solution based on a Hankel type operator, which parallels the Hankel operator known in the $H^\infty$ control theory.
%
\section{TV Allpass Property of the Optimum}
\label{allpass}
The dual space of $\Eu{B}(\ell^2, \; \ell^2)$, $\Eu{B}(\ell^2, \; \ell^2)^\star$, is given by the space, $\Eu{B}(\ell^2, \; \ell^2)^\star \simeq {\cal C}_1 \oplus_1
{\cal K}^\perp$,
where ${\cal K}^\perp$ is the annihilator of ${\cal K}$, the space of compact operators on $\ell^2$, and the symbol
$\oplus_1$ means that if $\Phi \in {\cal C}_1 \oplus_1 {\cal K}^\perp$
then $\Phi$ has a unique decomposition as follows \cite{schatten}
\beqa
\Phi = \Phi_o + \Phi_T,  \;\;\; 
\| \Phi \| = \| \Phi_o \| + \| \Phi_T\|
\label{18}
\eeqa
where $\Phi_o \in {\cal K}^\perp$, and $\Phi_T$ is induced by the operator
$T \in {\cal C}_1$, i.e, $\Phi_T(A) = tr(AT)$, $A\in \Eu{B}(\ell^2, \; \ell^2)$. Banach space duality states that \cite{luen}
\beqa
\inf_{y \in J} \| x-y \| = \max_{\Phi \in J^\perp,
\; \|\Phi\| \leq 1} | \Phi (x) |
\label{19}
\eeqa
In our case $J = {\Eu B}_c(\ell^2, \; \ell^2)$. Since $\Eu{B}(\ell^2, \ell^2)^\star$
contains ${\cal C}_1$ as a subspace, then $\Eu{B}_c(\ell^2, \; \ell^2)^\perp$ contains
the preannihilator $S$, i.e., the following expression for the annihilator of $J$, $J^\perp$, is deduced
\beqa
J^\perp := \Eu{B}_c(\ell^2, \; \ell^2)^\perp = S \oplus_1 \Bigl({\cal K} \cap
\Eu{B}_c(\ell^2, \; \ell^2) \Bigl)^\perp \label{20}
\eeqa
A result in \cite{dav} asserts that if a linear functional $\Phi$ belongs to
the annihilator $J^\perp$ and $\Phi$ decomposes as $\Phi = \Phi_o + \Phi_T$, where $\Phi_o \in \Bigl({\cal K} \cap \Eu{B}_c(\ell^2, \; \ell^2) \Bigl)^\perp$ and $\Phi_T \in S$, then both $\Phi_o \in \Eu{B}_c(\ell^2, \; \ell^2)^\perp\;\;\;
{\rm and}\;\;\; \Phi_T \in \Eu{B}_c(\ell^2, \; \ell^2)^\perp$ as well. We have then the
following \cite{cdc04_tv}
\beqa
\min_{Q \in {\Eu B}_c(\ell^2, \; \ell^2)} \| T_{2i}^\star T_1 T_{3ci}^\star - Q\| =
\max_{\br{c} \Phi_o \in ({\cal K} \cap {\Eu B}_c)^\perp,\; T \in S \\
\| \Phi_o\| + \| T\|_1 \leq 1 \er}
| \Phi_o (T_{2i}^\star T_1 T_{3ci}^\star) + tr(TT_{2i}^\star T_1 T_{3ci}^\star) |
\label{21}
\eeqa
If $\Phi_{opt} = \Phi_{opt,o} + \Phi_{T_{opt}}$ achieves the maximum in the RHS of
(\ref{21}), and $Q_o$ the minimum on the LHS, then the alignment condition in the dual
is immediately deduced from (\ref{21})
\beqa
|\Phi_{opt,o} (T_{2i}^\star T_1 T_{3ci}^\star) + tr(T_{opt}T_{2i}^\star T_1
T_{3ci}^\star) | =
\|T_{2i}^\star T_1 T_{3ci}^\star - Q_o  \| \; (\|\Phi_{opt,o}\|+\|T_{opt}\|_1)
\label{22}
\eeqa
Assume further that {\bf (A2)}: $T_{2i}^\star T_1 T_{3ci}^\star$ is a compact
operator. This is the case, for example, if $T_1$ is compact, then
$\Phi_{opt,o} (T_{2i}^\star T_1 T_{3ci}^\star) =0$
and the maximum in (\ref{21}) is achieved on $S$, that is the supremum in (\ref{16})
becomes a maximum. Compact operators include, for example, systems with impulse responses that have measurable and square integrable kernels.
\\
\\
It is instructive to note that in the LTI case assumption (A2) is the analogue of the assumption that $T_{2i}^\star T_1
T_{3ci}^\star$ is the sum of two parts, one part continuous on the unit circle and the other in $H^\infty$, in which case the optimum is allpass \cite{fran,zo}.
We would like to find the allpass equivalent for the optimum in the linear time varying case. This may be formulated by noting that flatness or
allpass condition in the LTI case means that the modulus of the optimum $|(T_{2i}^\star T_1 T_{3ci}^\star - Q_o) (\eio)|$ is constant at almost
all frequencies (equal to $\mu$). In terms of operator theory, note that the optimum viewed as a multiplication operator acting on $L^2$ or $H^2$,
changes the norm of any function in $L^2$ or $H^2$ by multiplying it by a constant (=$\mu$). 
That is, the operator achieves its norm at every $f \in L^2$ of unit $L^2$-norm.
This interpretation is generalized to the LTV case in the following Theorem, which part
of it first appeared in \cite{cdc04_tv} without a proof. In fact, the optimal cost
(\ref{4}) is shown to be a {\it partial isometry}.
\begin{theorem} \label{t2}
Under assumptions (A1) and (A2) there exists at least one optimal linear time varying
$Q_o \in {\Eu B}_c(\ell^2, \; \ell^2)$ that satisfies the following
\bi
\item[i)] the duality expression
\beqa
\mu = \|T_{2i}^\star T_1 T_{3ci}^\star -Q_o  \|
|tr (T_oT_{2i}^\star T_1 T_{3ci}^\star)| 
= \max_{n} \| (I-\widetilde{Q}_n) T_{2i}^\star T_1 T_{3ci}^\star \widetilde{Q}_n \|
\label{2400}
\eeqa
holds, where $T_o$ is some operator in $S$ and $\| T_o \|_1 =1$.
\item[ii)] and if
$\mu> \mu_{oo}:= \inf_{Q \in {\cal K}}
\|T_{2i}^\star T_1 T_{3ci}^\star -Q\|$, i.e. when the causality constraint is removed,
the following operator
\beqa
\frac{T_{2i}^\star T_1 T_{3ci}^\star -Q_o}{\|T_{2i}^\star T_1 T_{3ci}^\star -Q_o\|}
\label{2300}
\eeqa
is a {\it partial isometry}. That is, the optimum is an isometry on the range space of the operator $T_o$ in
$i)$. This is the time-varying counterpart of flatness. Note the condition $\mu > \mu_{oo}$ is sharp, in the sense that if it does not hold there exist $T_{2i}^\star T_1 T_{3ci}^\star$ and $Q_o$ such that $\frac{T_{2i}^\star T_1 T_{3ci}^\star -Q_o}{\|T_{2i}^\star T_1 T_{3ci}^\star -Q_o\|}$
is not a partial isometry.
\ei
\end{theorem}
{\bf Proof.} i) Identity (\ref{2400}) is implied by the previous argument that the supremum in
(\ref{16}) is achieved by some $T_o$ in $S$ with trace-class norm equal to 1. Combining
this result with Corollary 16.8 in \cite{dav} (see also \cite{arveson}),
which asserts that
\beqa
\|T_{2i}^\star T_1 T_{3ci}^\star -Q_o  \| = \sup_n \|(I-\widetilde{Q}_n) T_{2i}^\star T_1 T_{3ci}^\star \widetilde{Q}_n \|
\eeqa
shows in fact that the supremum w.r.t. $n$ is achieved proving that (\ref{2400}) holds.
\\
ii) The operator $T_o^\star T_o$ is self-adjoint and compact, and therefore admits a spectral representation
$T_o^\star T_o = \sum_j \lambda_j \phi_j \otimes \phi_j$, where $\{\phi_j\}$ form an orthonormal basis for
$\ell^2$ consisting of eigenvectors of $T_o^\star T_o$, and $\lambda_j$ are its necessarily positive real
eigenvalues \cite{schatten}. It follows that
\beqa
\sum_j <T_o^\star T_o \phi_j, \phi_j > = tr(T_o^\star T_o)^{\frac{1}{2}}
=\sum_j \lambda_j= \| T_o \|_1
\eeqa
where we used the fact that $T_o^\star T_o \phi_j = \lambda_j \phi_j \neq 0$.
Let the polar decomposition of $T_o$ be $T_o = U(T_o^\star T_o)^{\frac{1}{2}}$,
where $U$ is an isometry on the set $\{\phi_j\}$, and $(T_o^\star T_o)^{\frac{1}{2}}$ the "square root"
of $T_o^\star T_o$. Now $\sqrt{\lambda_j}$ being the non-zero singular values of $T_o$, are also the
eigenvectors of $(T_o^\star T_o)^{\frac{1}{2}}$. It follows that for $\lambda_j \neq 0$, we have
\beqa
 \frac{1}{\sqrt\lambda_j} T_o \phi_j = \frac{1}{\sqrt{\lambda_j}}
U (T_o^\star T_o)^\frac{1}{2}\phi_j = \frac{1}{\sqrt{\lambda_j}} U {\sqrt{\lambda_j}} \phi_j
=  U \phi_j
\eeqa
so $\{U\phi_j\}$ is an orthonormal set which spans the range of $T_o$. Call
$\psi_j := U \phi_j$, then $T_o$ can be written as $T_o = \sum_j \lambda_j \phi_j\otimes \psi_j$
and
\beqa
\| T_o\|_1 &=& \sum_j \lambda_j =1
\\
\mu &=& \left|tr \Bigl((T_{2i}^\star T_1 T_{3ci}^\star - Q_o) T_o \Bigl)\right|
=\left| \sum_j \lambda_j <\phi_j, (T_{2i}^\star T_1 T_{3ci}^\star - Q_o)\psi_j >
\right| \;\; {\mathrm{by \; definition\; of \; the\;trace}} \nonumber  \\
&\leq& \sum_j \lambda_j \|\phi_j\|_2
\|(T_{2i}^\star T_1 T_{3ci}^\star - Q_o)\psi_j\|_2
\leq \sum_j \lambda_j  \|T_{2i}^\star T_1 T_{3ci}^\star - Q_o\| \|\psi_j\|_2 \label{cauchy} \\
&\leq& \sum_j \lambda_j   \|T_{2i}^\star T_1 T_{3ci}^\star -Q_o  \|
\leq \| T_o\|_1 \|T_{2i}^\star T_1 T_{3ci}^\star - Q_o\|  = \mu \label{orth}
\eeqa
where the first inequality in (\ref{cauchy}) follows from the Cauchy-Schwartz inequality, the second from a standard property of
the induced norm, and the fact that $\{\phi_j\}$ and $\{\psi_j\}$ are orthonormal sets. Hence equality must hold throughout yielding
\beqa
\|(T_{2i}^\star T_1 T_{3ci}^\star - Q_o)\psi_j\|_2  =
\|T_{2i}^\star T_1 T_{3ci}^\star - Q_o\| \|\psi_j\|_2
\eeqa
for each $\psi_j$, that is,
$T_{2i}^\star T_1 T_{3ci}^\star - Q_o$ attains its norm on
each $\psi_j$, it must then attain its norm everywhere on the span of $\{\psi_j\}$, so that
$\frac{T_{2i}^\star T_1 T_{3ci}^\star - Q_o}{\|T_{2i}^\star T_1 T_{3ci}^\star - Q_o\|}$
is {\it an isometry on the range of} $T_o$.
\\
Note that if condition $\mu > \mu_{oo}$ does not hold it is straightforward to find
a counter example such that the optimal is not a partital isometry. Choose for example
$T_{2i}^\star T_1 T_{3ci}^\star$ to be a strict upper triangular which is not a partial isometry, then an optimal $Q_o$ is the zero matrix. The optimum is equal to  $\|T_{2i}^\star T_1 T_{3ci}^\star \|$ which does not correspond to a partial isometry.  $\blacksquare$
\\
\\
Identity (\ref{2300}) represents the allpass condition in the time-varying case, since it corresponds
to the allpass or flatness condition in the time-invariant case for the standard optimal $H^\infty$
problem. In the next section, we show that the search over $Q$ can be restricted to compact
operators. This is achieved by introducing a Banach space notion known as M-ideals \cite{dav}.
\section{M-Ideals and Compact Youla TV Parameters}
\label{mideal}
In this section we show that the optimal TV $Q$ is compact. This will allow approximations by finite rank operators to any desired accuracy, and consequently bring in finite dimensional convex optimizations to carry out the computations. In order to achieve our objective we rely on the concept of $M$-ideals which is discussed next.
\\
Following \cite{dav} we say that a closed subspace $M$ of a Banach space $B$ is an
$M$-ideal if there exists a linear projection $\Pi: \; B^\star \longrightarrow M^\perp$ from
the dual space of $B$, $B^\star$ onto the annihilator of $M$, $M^\perp$ in $B^\star$,
such that for all $b^\star \in B^\star$, we have
\beqa
\| b^\star \| = \| \Pi b^\star \| + \|b^\star - \Pi b^\star \|
\label{25}
\eeqa
In this case $M^\perp$ is called an $L$-summand of $B^\star$. The range $N$ of $(I-\Pi)$
is a complementary subspace of $M^\perp$, and $B^\star = M^\perp \oplus_1 N$. A basic
property of $M$-ideals is that they are {\it proximinal}, that is,
for every $b \in B$, there is an $m_o$ in $M$ such that $
\inf_{m \in M} \| b-m\| = \|b-m_o\|$.
Under assumption (A2) we generalize Lemma 1.6. in \cite{garnett} to causal LTV systems, i.e., the space $\Eu B_c(\ell^2, \ell^2)$. Recall
that Lemma 1.6 states that if $f$ is a function continuous on the unit circle, i.e.,
$f \in C$, then
\beqa
\inf_{g \in H^\infty} \| f-g\|_{\infty} = \inf_{g \in A} \|f -g \|_{\infty}
\label{26b}
\eeqa
where $A$ is the disk algebra, i.e., the space of analytic and continuous function on
the unit disk $A = H^\infty \cap C$. That is it suffices to restrict the search to
functions in $A$.
To generalize (\ref{26b}) to causal LTV systems put $B :=\cal K$,
$M := \Eu B_c(\ell^2, \ell^2)$, and show that for $b \in \cal K$ we have
\beqa
\inf_{m \in \Eu B_c(\ell^2, \ell^2)} \| b-m\| =
\inf_{m \in \Eu B_c(\ell^2, \ell^2) \cap \cal K} \|b - m\|
\label{2000}
\eeqa
By Theorem 3.11 in \cite{dav}, $\Eu B_c(\ell^2, \ell^2) \cap \cal K$ is weak$^\star$
dense in $\Eu B_c(\ell^2, \ell^2)$. By Theorem 11.6 and Corollary 11.7 in \cite{dav},
$\Eu{B}_c(\ell^2, \ell^2) \cap {\cal{K}}$, is an $M$-ideal in $\Eu{B}_c(\ell^2, \ell^2)$,
and the quotient map $q_1: \; \Eu{B}_c(\ell^2, \ell^2) /\bigl(\Eu{B}_c(\ell^2, \ell^2)
\cap {\cal{K}}\bigl) \longmapsto {\Eu{B}}_{c}(\ell^{2}, \ell^{2})+{\cal{K}}/
{\Eu{B}}_{c}(\ell^2, \ell^2)$
%
%
is isometric. Likewise, the quotient map $q_2: \; {\cal{K}} /\bigl({\Eu{B}}_{c}(\ell^{2},
\ell^{2})\cap{\cal{K}}\bigl)\longmapsto {\Eu{B}}_{c}(\ell^2, \ell^2) + {\cal{K}} /
{\Eu{B}}_{c}(\ell^{2}, \ell^{2})$, is isometric. And the identity (\ref{2000}) holds.
In our case under assumption
{\bf (A2)} $b =T_{2i}^\star T_1 T_{3ci}^\star \in \cal K$, and $m=Q$
yields
\beqa
\inf_{Q \in \Eu B_c(\ell^2, \ell^2)}\| T_{2i}^\star T_1 T_{3ci}^\star-Q\|
\nonumber 
=
\inf_{Q \in \Eu B_c(\ell^2, \ell^2) \cap \cal K} \|T_{2i}^\star T_1 T_{3ci}^\star - Q\|
=
 \|T_{2i}^\star T_1 T_{3ci}^\star - Q_o\|
\label{28}
\eeqa
for some optimal $Q_o \in \Eu B_c(\ell^2, \ell^2) \cap \cal K$. That is under (A2)
the optimal $Q$ is compact, and thus it suffices to restrict the search in (\ref{28})
to causal and compact parameters $Q$.
\\
Now we turn our attention to studying questions regarding the uniqueness of the optimal
solution.
\section{On the Uniqueness of the Optimal TV Controller}
\label{uniqueness}
In this section under assumptions (A1), (A2), and $\mu > \mu_{oo}$, we prove that the optimal TV controller is essentially unique.
Since there is an one-to-one onto correspondence between optimal controllers and Youla parameters, it suffices to show that the latter is essentially unique.
Suppose by way of contradiction that there are two optimal compact Youla parameters $Q_1$ and $Q_2$ such that
\beqa
\mu = \|T_{2i}^\star T_1 T_{3ci}^\star -Q_1  \| =   \|T_{2i}^\star T_1 T_{3ci}^\star -Q_2  \|
\eeqa
and $T_{2i}^\star T_1 T_{3ci}^\star -Q_1$ and $T_{2i}^\star T_1 T_{3ci}^\star -Q_1$ are both isometries on the range of $T_o$, where $T_o$ is
the optimal dual operator in Theorem \ref{t2}. In this case, any convex combination of
$Q_1$ and $Q_2$
is also optimal since $
\|  T_{2i}^\star T_1 T_{3ci}^\star + \lambda Q_1 +(1-\lambda) Q_2 \| \geq \mu,\;\;
\forall \lambda \in (0, \; 1)$, and
\beqa
&&\|  \lambda (T_{2i}^\star T_1 T_{3ci}^\star)
+(1-\lambda) (T_{2i}^\star T_1 T_{3ci}^\star) +
\lambda Q_1 +(1-\lambda) Q_2 \| \nonumber \\
&\leq&  \lambda \|T_{2i}^\star T_1 T_{3ci}^\star -Q_1  \| +
(1-\lambda) \|T_{2i}^\star T_1 T_{3ci}^\star -Q_2  \|  =\mu
\eeqa
implying that
\beqa
\|  T_{2i}^\star T_1 T_{3ci}^\star + \lambda Q_1 +(1-\lambda) Q_2 \| = \mu,
\;\; \forall \lambda \in [0, \; 1]
\eeqa
In particular for $Q^\star = \frac{Q_1+Q_2}{2}$, let
$\Gamma:=T_{2i}^\star T_1 T_{3ci}^\star - Q^\star$, and note that
\beqa
\| \Gamma \| = \left\| \Gamma \pm \frac{Q_1 -Q_2}{2} \right\| = \mu
\eeqa
By the parallelogram law we have for all $x \in \ell^2$,
\beqa
\left\| \left(\Gamma - \frac{Q_1-Q_2}{2}\right) T_o x \right\|_2^2 +
\left\| \left(\Gamma + \frac{Q_1-Q_2}{2}\right) T_ox \right\|_2^2 =
2\|\Gamma (T_o x)\|_2^2 +
2 \left\| \frac{Q_1-Q_2}{2} (T_ox) \right\|_2^2
\eeqa
Since each of the operators $\Gamma$ and $\Gamma \pm \frac{Q_1-Q_2}{2}$ is an isometry on the range of $T_o$, we have
\beqn
\left\| \left(\Gamma \pm \frac{Q_1-Q_2}{2}\right) T_o x \right\|_2 = \|\Gamma T_ox\|_2 =\|T_ox\|_2
\eeqn
which implies that $\left\| \frac{Q_1-Q_2}{2} (T_ox) \right\|=0$, $\forall x\in \ell^2$, that is, we have necessarily $Q_1 \equiv Q_2$, showing uniqueness on the range of $T_o$. $\blacksquare$
\\
Next, we show that the operator $T_o$ in the dual maximization is unique. To see this, again suppose that there exists another operator $T_o^\prime$ such that
\beqa
\|T_{2i}^\star T_1 T_{3ci}^\star -Q_o\| = tr\Bigl((T_{2i}^\star T_1 T_{3ci}^\star -Q_o)T_o\Bigl)
= tr\Bigl((T_{2i}^\star T_1 T_{3ci}^\star -Q_o)T_o^\prime\Bigl)
\eeqa
By the proof of Theorem \ref{t2} there exist orthonormal sequences
$\{\phi_i\}$ and $\{\psi_j\}$ such that
\beqa
T_o = \sum_j \lambda_j \phi_j \otimes \psi_j, \;\;\;\;
\| T_o\|_1 = \sum_j \lambda_j =1.
\eeqa
Moreover since $
(T_{2i}^\star T_1 T_{3ci}^\star -Q_o)T_o = \sum_j \lambda_j \psi\otimes
(T_{2i}^\star T_1 T_{3ci}^\star -Q_o) \psi_j$,
we have then $
tr\Bigl((T_{2i}^\star T_1 T_{3ci}^\star -Q_o)T_o\Bigl)= \sum_j \lambda_j <\phi_j,
(T_{2i}^\star T_1 T_{3ci}^\star -Q_o)\psi_j>$, and $\|T_{2i}^\star T_1 T_{3ci}^\star -Q_o\|= \mu$ yields $ <\phi_j, (T_{2i}^\star T_1 T_{3ci}^\star -Q_o)\psi_j> = \mu,\;\;\; \forall j$. That is, $
(T_{2i}^\star T_1 T_{3ci}^\star -Q_o)\psi_j = \mu \phi_j, \;\;\forall j$. The latter shows that $
(T_{2i}^\star T_1 T_{3ci}^\star -Q_o)T_o = \bigl(T_o^\star T_o\bigl)^{\frac{1}{2}}$.
Similarly, $(T_{2i}^\star T_1 T_{3ci}^\star -Q_o)T_o^\prime = \bigl(T_o^{\prime\star}
T_o^\prime\bigl)^{\frac{1}{2}}$.
\\
Now suppose that $R:=T_o-T_o^\prime \neq0$, $R\in S$, i.e. $R$ is strictly causal,
then $\| T_o\|_1 = \|R+T_o^\prime\|_1 =\|T_o^\prime\|_1=1$, and
$tr\bigl((T_{2i}^\star T_1 T_{3ci}^\star -Q_o)(T_o^\prime +R)\bigl) = tr\bigl((T_{2i}^\star T_1 T_{3ci}^\star -Q_o)T_o^\prime\bigl)$.
Therefore, $\|T_o^\prime\|_1 \mu = \|T_o^\prime +R\|_1\mu = tr\bigl((T_{2i}^\star T_1 T_{3ci}^\star -Q_o)T_o^\prime\bigl)$, and then
\beqa
(T_{2i}^\star T_1 T_{3ci}^\star -Q_o) R &=& (T_{2i}^\star T_1 T_{3ci}^\star -Q_o) (T_o^\prime+R) -
(T_{2i}^\star T_1 T_{3ci}^\star -Q_o)T_o^\prime
\\
&=& \bigl((T_o^\prime +R)^\star (T_o^\prime +R)\bigl)^{\frac{1}{2}} -
\bigl(T_o^\prime T_o^\prime\bigl)^{\frac{1}{2}}
\label{28a}
\eeqa
showing that $(T_{2i}^\star T_1 T_{3ci}^\star -Q_o) R$ is self-adjoint. But
$(T_{2i}^\star T_1 T_{3ci}^\star -Q_o) R$ is strictly-causal, that is, belongs to $S$,
likewise for its adjoint $\bigl((T_{2i}^\star T_1 T_{3ci}^\star -Q_o) R\big)^\star \in S$. This implies that
\beqa
\bigl((T_{2i}^\star T_1 T_{3ci}^\star -Q_o) R\big)^\star \bigl((T_{2i}^\star T_1 T_{3ci}^\star -Q_o) R\big) \in S
\eeqa
Since $\EuBc$ is the annihilator of $S$ it follows that
\beqa
tr\Bigl(\bigl((T_{2i}^\star T_1 T_{3ci}^\star -Q_o) R\big)^\star
\bigl((T_{2i}^\star T_1 T_{3ci}^\star -Q_o) R\big)\Bigl) =0
\eeqa
Thus, $(T_{2i}^\star T_1 T_{3ci}^\star -Q_o) R =0$, and from (\ref{28a}) we have $
\bigl((T_o^\prime +R)^\star (T_o^\prime +R)\bigl)^{\frac{1}{2}} = \bigl(T_o^\prime T_o^\prime\bigl)^{\frac{1}{2}}$. Since $T_{2i}^\star T_1 T_{3ci}^\star -Q_o$ is an isometry on the ranges of $T_o^\prime$ and $T_o=T_o^\prime +R$, we have for every $x\in \ell^2$,
\beqa
\|(T_o^\prime +R)x\|_2 &=& \|(T_{2i}^\star T_1 T_{3ci}^\star -Q_o)( T_o^\prime +R)x \|_2 \\
&=& \|(T_{2i}^\star T_1 T_{3ci}^\star -Q_o)T_o^\prime x\|_2 = \|T_o^\prime x\|_2
\eeqa
In the same manner, we get $\|(T_o^\prime -R)x \|_2 = \|T_o^\prime x\|_2$. Next, apply the parallelogram law to get the equality $
 \|(T_o^\prime +R)x \|_2^2 +   \|(T_o^\prime -R)x \|_2^2 = 2  \|T_o^\prime x \|_2^2 +2\|Rx \|_2^2$, from which we deduce that$\|Rx \|_2 =0$ for all $x\in \ell^2$, that is, we must have $R \equiv 0$ and then $T_o \equiv T_o^\prime$, showing uniqueness of $T_o$. $\blacksquare$
\\
\\
In the next section, we relate our problem to an LTV operator analogous to the Hankel operator, which is known to solve the standard
optimal $H^\infty$ control problem \cite{fran,zhou}. 
\section{Triangular Projections and Hankel Forms}
\label{hankel}
Let ${\cal C}_2$ denote the class of compact operators on $\ell^2$ called the
Hilbert-Schmidt or Schatten 2-class \cite{schatten,dav} under the norm,
\beqa
\| A \|_2 \; :=\; \Bigl(tr (A^\star A)\Bigl)^{\frac{1}{2}},\;\; A \in {\cal C}_2
\eeqa
Define the space ${\cal A}_2 := {\cal C}_2 \cap {\Eu{B}}_{c}({\ell^{2}}, {\ell^{2}})$,
then ${\cal A}_2$ is the space of causal Hilbert-Schmidt operators. The orthogonal projection ${\cal P}$ of ${\cal C}_2$ onto ${\cal A}_2$ is the lower triangular truncation, and is analogous to the standard positive Riesz projection
(for functions on the unit circle). Following \cite{power} an operator $X$ in ${\cal B}(\ell^2, \ell^2)$ determines a Hankel
operator $H_X$ on ${\cal A}_2$ if
\beqa
H_X A = (I-{\cal P}) X A, \;\;\; {\rm for}\; A \in {\cal A}_2
\eeqa
In the sequel we show that $\mu$ is equal to the norm of a particular LTV Hankel operator. To this end, we need first to characterize all atoms, denoted $\Delta_n$, of ${\Eu B_c}(\ell^2, \ell^2)$ as $\Delta_n :=\widetilde{Q}_{n+1}-\widetilde{Q}_n,\;\;\; n=0,1,2, \cdots$. Write ${\cal C}^+$ for the set of operators $A\in {\Eu B_c}(\ell^2, \ell^2)$ for which $\Delta_n A \Delta_n =0,\;\; n=0,1,2, \cdots$
and let ${\cal C}_2^+ := {\cal C}_2 \cap \cal {\cal C}^+$.
In \cite{power,power1} it is shown that any operator $A$ in ${\Eu{B}}_c(\ell^2, \ell^2)
\cap {\cal C}_1$ admits a {\it Riesz factorization}, that is, there exist operators
$A_1$ and $A_2$ in ${\cal A}_2$ such that $A$ factorizes as
\beqa
A = A_1 \; A_2, \;\;\; {\rm and} \;\; \| A \|_1 &=& \|A_1\|_2\; \|A_2 \|_2 \label{27}
\eeqa
%
A Hankel form $[\cdot \; , \; \cdot]_B$ associated to a bounded linear operator
$B \in {\Eu B}(\ell^2,\; \ell^2)$ is defined by \cite{power,power1}
\beqa
[A_1,\; A_2]_B = tr(A_1 B A_2),\;\; A_1,\; A_2 \in {\cal A}_2  \label{27a}
\eeqa
Since any operator in the preannihilator $S$ belongs also to
${\Eu{B}}_c(\ell^2, \ell^2) \cap {\cal C}_1$, then any $A \in S$ factorizes as in
(\ref{27}). And if $\| A\|_1 \leq 1$, as on the RHS of (\ref{16}), $A\in S$, the operators
$A_1$ and $A_2$ both in ${\cal A}_2$ can be chosen such that $\| A_1\|_2 \leq 1$, $\|A_2\|_2 \leq 1$, and $\Delta_n A_1\Delta_n =0, \;\; n=0,1,2, \cdots$, that is, $A_1 \in {\cal C}_2^+$ \cite{power1}.
\\
Now write ${\cal P}_+$ for the orthogonal projection with range ${\cal A}_2 \cap {\cal C}^+$. Introducing the notation $(B_1, B_2) = tr(B_2^\star \; B_1)$, and
the Hankel form associated to $B:=T_{2i}^\star T_1 T_{3ci}^\star$, we have by a result in \cite{power1},
\beqa
[A_1, A_2]_{B} = tr(B A_2 A_1) = (A_1, (B A_2)^\star )
= ({\cal P}_+ \; A_1, (BA_2)^\star ) =(A_1, {\cal P}_+ (BA_2)^\star ) =
(A_1, H_B^\star A_2)
\eeqa
 where $H_B$ is the Hankel operator $(I-{\cal P})B{\cal P}$ associated with
$B:=T_{2i}^\star T_1 T_{3ci}^\star$. The Hankel operator $H_B$ belongs to
the Banach space of bounded linear operators on ${\cal C}_2$. Furthermore, we have
\beqa
\| H_{T_{2i}^\star T_1 T_{3ci}^\star}\| = \sup_{\|A_2\|_2 \leq 1, A_2\in {\cal C}_2} \| H_B A_2\|_2= \sup_{\br {c}\|A_2\|_2 \leq 1, A_2\in {\cal C}_2 \\ \| A_1\|_2 \leq 1, A_1 \in {\cal C}_2^+ \er} (A_1, H_B^\star A_2) \label{150}
\eeqa
We have then the following Theorem which relates the optimal performance $\mu$ to the
induced norm of the Hankel operator $H_{T_{2i}^\star T_1 T_{3ci}^\star}$. The Theorem was
announced in \cite{cdc04_tv} without a proof.
\begin{theorem} \label{t3} Under assumptions (A1) and (A2) the following hold
\beqa
\mu = \| H_{T_{2i}^\star T_1 T_{3ci}^\star} \| = \|(I-{\cal P})T_{2i}^\star T_1 T_{3ci}^\star{\cal P} \|
\eeqa
\end{theorem}
{\bf Proof.} Since by the previous discussion any operator $T \in S$ can be
factored as $T = A_1 A_2$, where $A_1 \in {\cal C}_2^+$, $A_2 \in {\cal A}_2$,
$\| A_1 \|_2 = \|A_2\|_2 \leq 1$, and $\| T \|_1 = \|A_1 \|_2 \; \|A_2\|_2$,
the duality identity (\ref{16}) yields
\beqn
\mu &=&  \sup_{T \in S, \; \| T\|_1 \leq 1} |tr(TT_{2i}^\star T_1 T_{3ci}^\star)|
= \sup_{\br {c}\|A_2\|_2 \leq 1, A_2\in {\cal A}_2 \\ \| A_1\|_2 \leq 1, A_1
\in {\cal C}_2^+ \er} |tr( T_{2i}^\star T_1 T_{3ci}^\star A_2 A_1)|\\
&=& \| H_{T_{2i}^\star T_1 T_{3ci}^\star} \|, \;\;\;\;\;\;{\rm by}\;\; (\ref{150})
\eeqn
By Theorem 2.1. \cite{power1} the Hankel operator is a compact operator if and only if
$B$ belongs to the space ${\Eu B}_c(\ell^2, \ell^2)+ {\cal{K}}$. It follows in our
case that under assumption (A2) $H_{T_{2i}^\star T_1 T_{3ci}^\star}$ is a compact
operator on ${\cal A}_2$. A basic property of compact operators on Hilbert spaces is that they have maximizing vectors, that is, there exists an $A \in {\cal A}_2$, $\| A\|_2 =1$
such that $
\| H_{T_{2i}^\star T_1 T_{3ci}^\star}\| = \| H_{T_{2i}^\star T_1 T_{3ci}^\star} A \|_2$, that is, a $A$ achieves the norm of $H_{T_{2i}^\star T_1 T_{3ci}^\star}$. We can then
deduce from (\ref{16}) an expression for the optimal TV Youla parameter $Q_o$ as follows
\beqn
\|H_{T_{2i}^\star T_1 T_{3ci}^\star}\| = \| H_{T_{2i}^\star T_1 T_{3ci}^\star} A \|_2
= \| H_{T_{2i}^\star T_1 T_{3ci}^\star - Q_o} A \|_2 = \| \bigl(I-{\cal P}\bigl) \bigl(T_{2i}^\star T_1 T_{3ci}^\star A - Q_o A\bigl)
\|_2 \nonumber \\
\leq \|T_{2i}^\star T_1 T_{3ci}^\star A - Q_o A  \|_2
\leq \|T_{2i}^\star T_1 T_{3ci}^\star - Q_o  \|
\|A\|_2 \leq \| T_{2i}^\star T_1 T_{3ci}^\star - Q_o  \| = \| H_{T_{2i}^\star T_1 T_{3ci}^\star}\|
\eeqn
All terms must be equal, and then
\beqa
\|\bigl(I-{\cal P}\bigl)
\bigl(T_{2i}^\star T_1 T_{3ci}^\star A - Q_o A\bigl)
\|_2 = \|T_{2i}^\star T_1 T_{3ci}^\star A - Q_o A  \|_2
\eeqa
Since $
T_{2i}^\star T_1 T_{3ci}^\star A - Q_o A = \bigl(I-{\cal P}\bigl)
\bigl(T_{2i}^\star T_1 T_{3ci}^\star A - Q_o A\bigl) +  {\cal P}
\bigl(T_{2i}^\star T_1 T_{3ci}^\star A - Q_o A\bigl)$
and $ {\cal P} \bigl(T_{2i}^\star T_1 T_{3ci}^\star A - Q_o A\bigl) =0$. The TV optimal $Q_o$ can then be computed from the following operator identity
\beqa
Q_o A =  T_{2i}^\star T_1 T_{3ci}^\star A - H_{T_{2i}^\star T_1 T_{3ci}^\star}A
\eeqa
The upshot of these methods is that they lead to the computation of $\mu$
within desired tolerances by solving two finite dimensional convex optimizations.
\section{Numerical Computation of the Optimal Solution}
\label{numerical}
In this section a numerical solution based on duality theory is developed.
If $\{e_n \; :\; n=0,1,2, \cdots\}$ is the standard orthonormal basis in $\ell^2$,x
then $\widetilde{Q}_n \ell^2$ is the linear span
of $\{ e_k\;:\; k= n+1, n+2, \cdots\}$. The matrix representation of
$A \in {\Eu B}_c(\ell^2, \ell^2)$ w.r.t.
this basis is lower triangular. Note $P_n =I-\widetilde{Q}_n \longrightarrow I$ as $n \longrightarrow \infty$
in the strong operator topology (SOT). If we restrict the minimization in
(\ref{16}) over $Q \in {\Eu B}_c(\ell^2, \ell^2)$ to the span of $
\{e_n \; :\; n=0,1,2, \cdots, N\}$, that is, $P_N \ell^2 =: \ell^2_N$,
this yields a finite dimensional convex optimization problem in lower triangular matrices
$Q_N$ of dimension $N$, that is,
\beqa
Q_N = \left( \br{ccccc} Q_{11} & 0 & 0 & \cdots & \cdots \\
Q_{21}& Q_{22} & 0 & \cdots & \cdots\\
\vdots & \vdots & \ddots & \cdots &\cdots
\\
Q_{N1} & Q_{N2} & Q_{N3} & \cdots & Q_{NN}
\er
\right)
\eeqa
and the optimization
\beqa
\mu_N:=\inf_{Q_N \in \Eu{B}_c(P_N\ell^2, P_N\ell^2)}
\| T_{2i}^\star T_1 T_{3ci}^\star - Q_N  \|
\label{25b}
\eeqa
where
{\scriptsize
\beqa
T_{2i}^\star T_1 T_{3ci}^\star -Q_N= \left( \br{ccccccc}
T_{11} - Q_{11} & T_{12} & T_{13} & \cdots & T_{1N} & T_{1(N+1)} & \cdots
\\ T_{21}-Q_{21} & T_{22}- Q_{22} & T_{23} & \cdots & T_{2N} & T_{2(N+1)}& \cdots
\\ T_{31}-Q_{31} & T_{32}-Q_{32} & T_{33}-Q_{33} & \cdots & T_{3N} & T_{3(N+1)}& \cdots
\\ \vdots & \vdots & \vdots & \ddots & \vdots & \vdots & \cdots
\\ T_{N1}-Q_{N1} & T_{N2}-Q_{N2} & T_{N3} -Q_{N3}& \cdots & T_{NN}-Q_{NN} & T_{N(N+1)}& \cdots
\\
T_{(N+1)1} & T_{(N+1)2} & T_{(N+1)3} & \cdots & T_{(N+1)N} & T_{(N+1)(N+1)}& \cdots
\\
\vdots & \vdots & \vdots & \vdots & \vdots & \vdots & \cdots
\er
\right)
\eeqa
}
where $T_{ij}$ are fixed and correspond to the entries of $T_{2i}^\star T_1 T_{3ci}^\star$,
and $Q_{ij}$ are variable. The optimization $\mu_N$ overestimates $\mu$, but $\mu_N\downarrow\mu$, and
results in upper bounds and suboptimal TV parameters $Q_{N,o}$ and control laws, since $Q$ is restricted
to a proper subspace of ${\Eu B}_c(\ell^2, \ell^2)$. By Arveson distance formula \cite{arveson}
the minimization (\ref{25b}) is equal to
\beqa
\mu_N = \max_{1\leq n \leq N} \| (I-\widetilde{Q}_n)(T_{2i}^\star T_1 T_{3ci}^\star)\widetilde{Q}_n\|
\label{arv}
\eeqa
The degree of suboptimality can be computed
explicitly as follows:
Applying the same argument to the dual optimization on the RHS of
(\ref{16}), by restricting $S$ to the finite dimensional subspace
\beqa
S_N := \{T_N \in {\cal C}_1(P_N\ell^2,
P_N\ell^2):\; (I-\widetilde{Q}_n)T_N\widetilde{Q}_{n+1} =0,\;\; {\rm for\; all}\;\; 0\leq n \leq N-1\}
\label{25c}
\eeqa
In fact, with respect to the canonical basis $\{e_i\;,\; i=0, 1 ,\cdots\}$ for $P_N\ell^2$
the subspace $S_N$ is nothing but the space of $N\times N$ strictly lower triangular matrices.
\\
The dual optimization becomes
\beqa
\mu_{N}^\prime :=\sup_{T_N \in S_N, \; \| T_N\|_1 \leq 1}
\left| tr(T_N(T_{2i}^\star T_1 T_{3ci}^\star)
|_N\right|
\label{25d}
\eeqa
where $(T_{2i}^\star T_1 T_{3ci}^\star)|_N$ is the restriction of the operator
$T_{2i}^\star T_1 T_{3ci}^\star$ to $P_N\ell^2$. The supremum in (\ref{25d})
is in fact a maximum, since (\ref{25d}) is a maximization of a continuous linear functional over a compact
set. The optimization (\ref{25d}) is a finite variable constrained convex optimization in the (strictly
lowers entries of $T_N$), which yields lower bounds for $\mu$ since the dual optimization involves a supremum
rather than an infimum, i.e., $\mu_{N}^\prime \uparrow \mu$, and suboptimal $T_{N,o}$.
\\
The optimization problem (\ref{25d}) is in fact a semi-definite programming problem (SDP) since $\|T_N\|_1 \leq 1$ if and only
if there exist matrices $Y, Z \in \Bbb R^{N \times N}$ such that \cite{fazel,vander}
\beqa
\left(\br{cc} Y & T_N \\ T_N^T & Z \er \right) \geq 0, \;\;\;\; tr Y + tr Z \leq 2
\label{25x}
\eeqa
The optimization (\ref{25d}) becomes then
\beqa
&{\rm supremum \;\;}&tr(T_N(T_{2i}^\star T_1 T_{3ci}^\star) |_N) \label{25y} \\
&{\rm subject \;\; to\;\;}& (\ref{25x}) \nonumber \\
&& P_n T_N (I-P_{n+1}) = 0, \;\;\; n=0, 1, \cdots, N-1 \nonumber
\eeqa
where $P_n$ is the $N\times N$ truncation matrix having as its first $n$ columns the standards basis vectors $\{e_i\}$,
$i=1,2, \cdots, n$ and the remaining columns as $N\times 1$ zero vectors, that is,
\beqn
P_n =\left( e_1, e_2, \cdots, e_n, 0, \cdots, 0 \right)
\eeqn
Under assumption (A2) $T_{2i}^\star T_1 T_{3ci}^\star$ and the
optimal $Q_o$ are compact, and by a compactness argument we have that $Q_{N,o}\longrightarrow Q_o$
in the {\it operator topology}, i.e., $\| Q_o-Q_{N,o}\| \longrightarrow 0$ as $N \longrightarrow \infty$.
Likewise $T_{N,o} \longrightarrow T_o$ as $N\longrightarrow 0$ in the trace class topology. Since $P_N \longrightarrow I$ as $N \longrightarrow \infty$ in the SOT.
It is straightforward to show that these upper and lower bounds $\mu_N$ and $\mu_{N}^\prime$ converge to the optimum $\mu$ as $N \longrightarrow \infty$. These optimizations estimate $\mu$ within known tolerance and compute the corresponding LTV operators $Q_N$, which in turn allow the computation of LTV controllers $K$ through the Youla parametrization.
Solving such problems are then applications of finite variable convex programming techniques.
\\
Note the size of the truncations depends on the degree of accuracy desired. However, the size required for a particular accuracy can be exactly estimated by looking at the difference between the optimization (\ref{arv}) and the dual optimization (\ref{25d}). For {\it periodic systems} \cite{lall}, say of period $q$, it suffices to take the first $q$ vectors $\{e_n\}_{n=1}^q$ of the standard basis, i.e., $N=q$. In this case and for {\it finite horizon} problems
the finite dimensional convex optimizations yield {\it exactly} the optimal corresponding TV Youla parameter
$Q_o$ and hence the optimal TV controller.
\section{The Mixed Sensitivity Problem for LTV Systems}
\label{mixed}
The mixed sensitivity problem for stable plants \cite{dft,zhou} involves the sensitivity operator $T_{1} := \vect{W}{0}$, the complementary sensitivity operator $T_2=\vect{W}{V}P$ and $T_3:=I$ which are all assumed to belong to $\EuBct$, and is given by the optimization
\beqa
\mu_o = \inf_{Q \in \EuBc} \left\| \vect{W}{0}-\vect{W}{V}PQ \right\|
\label{151}
\eeqa
where $\| \cdot\|$ stands for the operator norm in $\EuBt$ given by
\beqa
\| B \| = \sup_{\|x\|_2 \leq 1, \; x \in \ell^2} \Bigl( \| B_1 x \|_2^2 +\| B_2 x\|_2^2 \Bigl)^{\frac{1}{2}},
\;\; B = \vect{B_1}{B_2}
\eeqa
The optimization problem (\ref{151}) can be expressed as a distance problem from the operator
$T_1$  to the subspace $\cm S = T_2 P \; \EuBc$ of $\EuBt$.  \\
To ensure closedness of $\cm S$, we assume that $W^\star W   +V^\star V > 0$, i.e.,
$W^\star W +V^\star V$ as an operator acting on $\ell^2$ is a positive operator.
Then there exists an outer spectral factorization $\Lambda_1 \in \EuBc$, invertible in
$\EuBc$ such that $\Lambda_1^\star  \Lambda_1 =  W^\star W  +V^\star V$ \cite{arveson,feintuch}.
Consequently, $\Lambda_1 P$ as a bounded linear operator in $\EuBc$ has an inner-outer factorization
$U_1 G$, where $U_1$ is inner and $G$ an outer operator defined on $\ell^2$ \cite{dav}.
\\
Next we assume ({\bf A3}) $G$ is invertible, so $U_1$ is unitary, and the operator $G$ and its inverse $G^{-1} \in \EuBc$. (A3) is satisfied when, for e.g., the outer factor of the plant is invertible. Let $R = T_2 \Lambda_1^{-1}U_1$, assumption (A3) implies that the operator $R^\star R \in \EuB$ has a bounded inverse, this ensures closedness of $\cm S$. According to Arveson (Corollary $2$, \cite{arveson}),
the self-adjoint operator $R^\star R$ has a spectral factorization of the form: $R^\star R =
\Lambda^\star \Lambda$, where $\Lambda,\; \Lambda^{-1}\in \; \EuBc$. Define $R_2 = R
\Lambda^{-1}$, then $R_2^\star R_2= I$, and  $\cm S$ has the equivalent representation, $\cm S = R_2 \EuBc$.
After "absorbing" $\Lambda$
into the free parameter $Q$, the optimization problem (\ref{151}) is then equivalent to:
\be
\mu_o = \inf_{Q \in \EuBc} \left\| T_1  - R_2Q \right\|
\label{152} 
\ee
To solve the TV optimization (\ref{152}) it suffices to apply the duality results of section
\ref{duality}. The latter yields the predual space of $\EuBt$ under trace duality, as the Banach
space isometrically isomorphic to
\beqa
\Eu C_1:= \left\{ B= \vect{B_1}{B_2}:
\; B_i \in {\cal{C}}_1\right\}
\eeqa
under the norm
\beqa
\|B\|^2_{11} := tr(B_1^\star B_1 + B_2^\star B_2)^{\frac{1}{2}}
\eeqa
The preannihilator $^\perp \cm S$ of $\cm S$ is characterized in the following Lemma.
\begin{lemma}
The preannihilator $^\perp \cm S$ can be computed as
\beqa
^\perp\cm S = R_2 S \oplus (I-R_2R_2^\star)
\label{152a}
\Eu C_1
\eeqa
where $\oplus$ denotes the direct sum of two subspaces.
\end{lemma}
{\bf Proof.}
To show (\ref{152a}) notice that $tr (\Phi^\star T) =0$, $\forall T$ in $^\perp\cm S$ for
$\Phi \in \EuBt$, is equivalent to $\Phi^\star (I-RR^\star) =0$ and $\Phi^\star R = A^\star$
for some $A\in B_c(\ell^2, \ell^2)$. This implies that  $\Phi^\star = A^\star R^\star$.
By taking the adjoints we get $\Phi=R_2A \in \cm S$.
\\
\\
The following lemma is a consequence of (\ref{15a}).
\begin{lemma}\label{l1} Under assumption (A3) there exists at least one optimal TV operator
$Q_o \in \EuBc$ s.t.
\beqa
\mu_o = \| T_1 - R_2 Q_o \| = \sup_{s \in ^\perp \cm S, \| s\|_{11} \leq 1} |trT_1^\star s |
\label{153}
\eeqa
\end{lemma}
Note that Lemma yields \ref{l1} not only shows existence of an optimal TV controller $K$ through the Youla
parameter $Q_o$  under assumption (A3), but also leads to two dual finite variable convex programming problems
using the same argument as in section \ref{numerical}, which solutions yield the optimal $Q_o$ within desired
tolerance. A TV allpass property holds also for the TV mixed sensitivity problem under the assumption that $T_1$ is
a compact operator in $\EuBct$ and $\mu_o > \mu_{oo}^\prime$, where $\mu_{oo}^\prime:= \inf_{Q\in \cal K} \| T_1 -R_2 Q \|$,
i.e., when the causality condition of $Q$ is removed. This generalizes the allpass property of section
\ref{allpass} and is summarized in the following lemma.
\begin{lemma}\label{l2}  If $T_1$ is a compact operator and $\mu_o > \mu_{oo}^\prime$ then the optimal mixed
sensitivity operator $\frac{T_1 -R_2 Q_o}{\mu_o}$ is a partial isometry
from $\ell^2$ onto $\ell^2 \times \ell^2$. The condition $\mu_o > \mu_{oo}^\prime$ is sharp,
in the sense that if it does not hold there exists $T_1$ and $R_2$ such that the optimum is not allpass.
\end{lemma}
{\bf Proof.} The proof is similar to the proof of Theorem 1 by working in the predual maximization.
The details are omitted. If $\mu_o > \mu_{oo}^\prime$ does not hold, take for example $V=1$ and $W=\frac{1-z}{2}$,
in which case $Q_o = 0$ showing that the optimum is not allpass.
\\
\\
Lemma \ref{l2} is the TV version of the same notion known to hold in the LTI case as shown
earlier in \cite{verma} using broadband theory and \cite{acc05_tv}.
\\
\\
Next we show that the mixed sensitivity optimization is equal to the norm of a certain TV Hankel-Topelitz operator. In order to do this we use a standard trick in \cite{feintuch}.
Since $R_2 = \vect{R_{21}}{R_{22}}$ is an isometry, so is $U:=\vect{R_2^\star}{I-R_2R_2^\star}$. Thus,
\beqn
 \mu_o =  \inf_{Q \in \EuBc} \left\| U \bigl( T_1 - R_2 Q\bigl) \right\|
 = \inf_{Q \in \EuBc} \left\| \vectt{R_{21}^\star W - Q}{(I-R_{21}R_{21}^\star )W}{-R_{22}
 R_{21}^\star W} \right\|
 \eeqn
Call $\Omega := \vect{(I-R_{21}R_{21}^\star )W}{-R_{22} R_{21}^\star W}$. By analogy with the
TV Hankel operator defined in section \ref{hankel}, we define the {\it TV Toeplitz operator}
$T_{\Omega^\star \Omega}$ associated to $\Omega^\star \Omega$ as:
\beqa
T_{\Omega^\star \Omega} A = {\cal P} \Omega^\star \Omega A,  \;\;\; {\rm for}\;\; A \in {\cal A}_2
\eeqa
The following Theorem solves the problem in terms of a TV mixed Hankel-Toeplitz operator.
\begin{theorem} \label{l3}  Under assumption (A3) the following holds
\beqa
\mu_o^2= \| H^\star_{R_{21}^\star W}  H_{R_{21}^\star W} + T_{\Omega^\star \Omega} \|
\label{hanktoe}
\eeqa
where $H_{R_{21^\star W}}$ is the Hankel operator associated to $R_{21}^\star W$, i.e., $H_{R_{21}^\star W} A= (I-{\cal P}) R_{21}^\star W A, \;\;\; {\rm for}\;\; A \in {\cal A}_2$.
\end{theorem}
{\bf Proof.} To prove (\ref{hanktoe}) define the following operator:
\beqa
\Gamma := \Pi_{{\cal{A}}_{2} \times {\cal{A}}_{2}  \ominus R_2 {\cal{A}}_{2}}  \vect{W}{0}
\eeqa
where $\Pi$ is the orthogonal projection from ${\cal{A}}_2 \times {\cal{A}}_2$ into the orthogonal
complement of ${\cal{A}}_{2}$, ${\cal{A}}_2 \times {\cal{A}}_2  \ominus R_2 {\cal{A}}_{2}$.
The orthogonal projection $\Pi$ can be shown to be given explicitly by \cite{djouadiSIAM}
\beqa
\Pi = I- R_2 {\cal{P}} R_2^\star
\label{proj}
\eeqa
It follows from the Commutant Lifting Theorem for nest algebra (Theorem 20.22 \cite{dav}) that
\beqa
\mu_o =   \inf_{Q \in \EuBc} \left\| T_1 - R_2 Q \right\| =\left\| \Gamma \right\|
\eeqa
Since $U$ is an isometry and thus preserves norms, operator pre-composition with $U$ and using the
explicit expression (\ref{proj}) of $\Pi$, straightforward computations show that $
\mu_o = \left\| \vect{H_{R_{21}^\star W}}{T_{\Omega}} \right\|$, and therefore (\ref{hanktoe}) holds. $\blacksquare$
\\
\\
Theorem \ref{l3} generalizes the solution of the mixed sensitivity problem in terms of a mixed Hankel-Toeplitz operator in the LTI case \cite{verma,ozbay,acc05_tv} to the TV case. This result also applies to solve the robustness problem of feedback systems in the gap metric \cite{GS2} in the TV case as outlined in \cite{feintuch,georgiou},
since the latter was shown in \cite{feintuch} to be equivalent to a special version of
the mixed sensitivity problem (\ref{151}).
\section{Conclusion}
\label{conclusion}
The optimal disturbance rejection problem for LTV systems involves solving a shortest distance minimization problem in the space of bounded linear operators. Dual and predual representations show existence of optimal TV controllers through the Youla parametrization. Under specific conditions, the optimal solution is compact and satisfies a``time-varying'' allpass or flatness condition. The proposed approach leads ``naturally'' to a numerical solution based on finite variable convex programming, which involves solving a SDP problem for the predual and a search over lower triangular matrices for the dual problem. Moreover, a solution based on a TV compact Hankel operator is proposed. The latter leads to an operator identity for the optimal TV Youla parameter $Q_o$. A generalization to the mixed-sensitivity problem is carried out that draws certain analogies with the LTI case.
\\
Future work includes investigation of the numerical solutions by semi-define programming and finding efficient algorithms to take advantage of the special structure (lower triangular) of $Q_o$ to solve numerically the dual problem.
\footnotesize

\end{document}